\documentclass[12pt]{amsart}%
\usepackage{amsmath}
\usepackage{amscd}
\usepackage{geometry}
\usepackage{amsfonts}
\usepackage{amssymb}
\usepackage{graphicx}%
\newtheorem{theorem}{Theorem}[section]
\theoremstyle{plain}

\newtheorem{lemma}[theorem]{Lemma}

\newtheorem{remark}{Remark}

\numberwithin{equation}{section}

\begin{document}
\title[On the dimension of $\mathcal{Z}$-sets]{On the dimension of $\mathcal{Z}$-sets}
\author{Craig R. Guilbault }
\address{Department of Mathematical Sciences, University of Wisconsin-Milwaukee,
Milwaukee, Wisconsin 53201}
\email{craigg@uwm.edu}
\author{Carrie J. Tirel}
\address{Department of Mathematical Sciences, University of Wisconsin-Fox Valley, 1478
Midway Road, Menasha, Wisconsin 54952}
\email{carrie.tirel@uwc.edu}
\thanks{This project was aided by a Simons Foundation Collaboration Grant awarded to
the first author.}
\date{July 13, 2013}

\begin{abstract}
We offer a short and elementary proof that, for a $\mathcal{Z}$-set $A$ in a
finite-dimensional ANR $Y$, $\dim A<\dim Y$. This result is relevant to the
study of group boundaries. The original proof by Bestvina and Mess relied on
cohomological dimension theory.

\end{abstract}
\maketitle

\section{Introduction\label{Section: Introduction}}

The notion of a $\mathcal{Z}$-set, developed by R.D Anderson \cite{An}, played
an important role in the infinite-dimensional topology boom of the 1970's (see
\cite{Ch}). More recently, work by Bestvina, Mess, Geoghegan and others (see
\cite{BM},\cite{Be},\cite{Ge}) has revived interest in $\mathcal{Z}$-sets. In
this more recent work, the point of view is the following: for an infinite
group $G$ acting nicely on an AR $X$, there is frequently a compactification
$\overline{X}=X\sqcup Z$, where $Z$ is a \textquotedblleft
boundary\textquotedblright\ for $G$ and a $\mathcal{Z}$-set in $\overline{X}$.
From there, well-known properties of $\mathcal{Z}$-sets reveal connections
between the group and its boundary. In contrast with earlier applications of
$\mathcal{Z}$-sets, the spaces involved in this new setting are usually
finite-dimensional and the precise dimensions are of utmost importance. For
example, \cite{BM} shows that, for torsion free hyperbolic $G$, the
topological dimension of the Gromov boundary $\partial G$ is always one less
than the (algebraically defined) cohomological dimension of $G$; the latter is
bounded above by the topological dimension of $X$. That result was expanded
upon in \cite{Be},\cite{GO}, and \cite{Dr}.

A useful lemma from \cite{BM} asserts that, in the above setting, the
dimension of $Z$ is strictly less than the dimension of $X$. The result was
obtained using \v{C}ech cohomology and the fact that, for finite-dimensional
spaces, the \textquotedblleft cohomological dimension\textquotedblright\ of a
space agrees with its topological dimension. The aim of this note is an
elementary proof of the same fact. No knowledge of cohomological dimension
theory is needed and no algebraic topology will be used. By \textquotedblleft
topological dimension\textquotedblright, we mean \textquotedblleft Lebesgue
covering dimension\textquotedblright, a concept to be defined shortly. Our aim
is the following:

\begin{theorem}
[{see \cite[Prop.2.6]{BM}}]\label{Main Theorem}If $A$ is a $\mathcal{Z}$-set
in a compact finite-dimensional metric ANR $Y$, then $\dim Y=\dim(Y-A)$ and
$\dim A<\dim Y$.
\end{theorem}

\begin{remark}
More general versions of this theorem are possible. For example, we make no
use of the ANR properties of $Y$ or $Y-A$; if $A$ is a subset of any compact
finite-dimensional metric space and it is possible to instantly homotope $Y$
off $A$ (see the first bullet point in
\S \ref{Section: Definitions and background}), then Theorem \ref{Main Theorem}
and its proof are still valid. One may also relax the compactness condition on
$Y$ to local compactness without much additional effort. We have chosen to
focus on the case of primary interest with the fewest technicalities.
\end{remark}

\section{Definitions and background\label{Section: Definitions and background}%
}

In this paper, all spaces are assumed to be separable metric. A locally
compact space $Y$ is an ANR (absolute neighborhood retract) if it can be
embedded as a closed subset of $%
\mathbb{R}
^{n}$ or $%
\mathbb{R}
^{\infty}$ so that it is a retract of one of its neighborhoods. A closed
subset $A$ of an ANR\ $Y$ is a $\mathcal{Z}$\emph{-set} if either of the
following equivalent conditions is satisfied:

\begin{itemize}
\item There exists a homotopy $H:Y\times\left[  0,1\right]  \rightarrow Y$
such that $H_{0}=\operatorname*{id}_{Y}$ and $H_{t}\left(  X\right)  \subseteq
Y-A$ for all $t>0$. (We say that $H$ \emph{instantly homotopes} $Y$ off from
$A$.)

\item For every open set $U$ in $Y$, $U-A\hookrightarrow U$ is a homotopy equivalence.
\end{itemize}

\noindent A $\mathcal{Z}$\emph{-com\-pact\-ific\-at\-ion} of a space $X$ is a
com\-pact\-ific\-at\-ion $\overline{X}=X\sqcup Z$ with the property that $Z$
is a $\mathcal{Z}$-set in $\overline{X}$. In this case, $Z$ is called a
$\mathcal{Z}$\emph{-boundary }for $X$. Implicit in this definition is the
requirement that $\overline{X}$ be an ANR. Since an open subset of an ANR is
an ANR, $X$ itself must be an ANR to be a candidate for $\mathcal{Z}%
$-com\-pact\-ific\-at\-ion. Hanner's Theorem \cite{Ha} ensures that every
com\-pact\-ific\-at\-ion $\overline{X}$ of an ANR $X$, for which $\overline
{X}-X$ satisfies either of the above bullet points is necessarily an ANR;
hence, it is a $\mathcal{Z}$-com\-pact\-ific\-at\-ion.

A collection $\mathcal{A}$ of subsets of a space $X$ has \emph{order} $k$ ($k$
an integer) if some $x\in X$ belongs to $k+1$ elements of $\mathcal{A}$ but
none belongs to more than $k+1$ elements of $\mathcal{A}$; it has \emph{mesh}
$\epsilon$ ($\epsilon\geq0$) if the diameters of its elements are bounded
above by $\varepsilon$. The collection $\mathcal{A}$ \emph{refines} a second
collection $\mathcal{B}$ of subsets of $X$ if, for each $A\in\mathcal{A}$
there exists $B\in\mathcal{B}$ such that $A\subseteq B$.

A space $X$ \emph{has dimension }$\leq k$ if, for every open cover
$\mathcal{U}$ of $X$, there exists an open cover $\mathcal{V}$ of $X$ that
refines $\mathcal{U}$ and has order $\leq k$. The smallest such integer (when
it exists) is called the \emph{dimension }of $X$ and is denoted $\dim X$. If
no such integer exists, $X$ is called \emph{infinite-dimensional}.

\begin{remark}
\emph{There are several equivalent definitions of dimension. We have chosen
the most elementary; it is frequently referred to as }Lebesgue covering dimension.
\end{remark}

The following is immediate from the existence of Lebesgue numbers.

\begin{lemma}
\emph{A compact metric space} $X$ \emph{has dimension} $\leq k$ \emph{if and
only if it admits open covers of order} $k$ \emph{with arbitrarily small
mesh.}
\end{lemma}

We call a map $f:X\rightarrow Y$ between metric spaces a $\delta$-$\epsilon
$-map if, for each $A\subseteq Y$ of diameter $\leq\delta$, $f^{-1}\left(
A\right)  $ has diameter $\leq\epsilon$.

\begin{lemma}
\label{delta-epsilon lemma} A compact metric space $X$ has dimension $\leq k$
if and only if, for every $\epsilon>0$, there exists $\delta>0$ and a $\delta
$-$\epsilon$-map from $X$ to a space $Y$ of dimension $\leq k$.

\begin{proof}
The identity map satisfies the forward implication, so we proceed to the
converse. Since every subspace of $Y$ also has dimension $\leq k$, we may
assume that $f$ is onto and $Y$ is compact. Choose an open cover $\mathcal{V}$
of $Y$ of order $\leq k$ and mesh $\leq\delta$. Then $\mathcal{U}=\left\{
f^{-1}\left(  V\right)  \mid V\in\mathcal{V}\right\}  $ has order $\leq k$ has
mesh $\leq\epsilon$.
\end{proof}
\end{lemma}

The final lemma of this section is less trivial than the previous two, but it
goes back to the beginnings of the subject. We refer the reader to the
original source \cite{Hu} for a proof.

\begin{lemma}
\label{Hurewicz lemma}For any nonempty locally metric compact space $X$,
$\dim\left(  X\times\left[  0,1\right]  \right)  =\dim X+1$.
\end{lemma}

\begin{remark}
\emph{Although elementary, Lemma \ref{Hurewicz lemma} is not immediate. There
exist compact spaces }$X$\emph{ and }$Y$\emph{ with }$\dim\left(  X\times
Y\right)  <\dim X+\dim Y$\emph{. The point of the Hurewicz paper is that this
does not happen when }$Y$\emph{ is 1-dimensional. As with our main theorem,
more general results are possible.}
\end{remark}

\section{Proof of Theorem \ref{Main Theorem}}

\begin{proof}
Suppose $A$ is a $\mathcal{Z}$-set in a compact finite-dimensional metric ANR
$Y$. Let $\epsilon>0$, and set $\delta:=\frac{\epsilon}{3}$. By definition
there is a homotopy $J:Y\times\lbrack0,1]\rightarrow Y$ satisfying
$J_{0}=\operatorname*{id}_{Y}$ and $J_{t}(Y)\subseteq Y-A$ for all $t>0$.

By compactness of $Y$, we may choose $T>0$ so that $d(y,J(y,t))<\frac
{\epsilon}{3}$ for all $y\in Y$ and $t\in\lbrack0,T]$. Define $H:Y\times
\lbrack0,1]\rightarrow Y$ by $H(y,t)=J(y,t\cdot T)$. Now $d(y,H(y,t))<\frac
{\epsilon}{3}$ for all $(y,t)\in Y\times\lbrack0,1]$, so $H_{1}$ is a $\delta
$-$\epsilon$-map from $Y$ to $Y-A$, and by Lemma \ref{delta-epsilon lemma}
$\dim Y=\dim(Y-A)$.

To show that $\dim A<\dim Y$, we will define $\delta^{\prime}>0$ and use the
homotopy $H$ to construct a $\delta^{\prime}$-$\epsilon$-map from
$A\times\lbrack0,1]$ into $Y$, where $A\times\lbrack0,1]$ is endowed with the
$\ell_{\infty}$-metric $d=\max\left\{  d_{1},d_{2}\right\}  $. From there, an
application of Lemma \ref{Hurewicz lemma} completes the proof.

Pick $k\in\mathbb{N}$ so that $\frac{1}{k}<\frac{\epsilon}{3}$. Then choose
open subsets $U_{0},U_{1},\ldots,U_{k+1}\subseteq Y$, and $t_{1},t_{2}%
,t_{3},\ldots,t_{k+1}\in\lbrack0,1]$ in the following way:\smallskip

\begin{itemize}
\item $U_{0}=\emptyset$ and $t_{1}=1$;\smallskip

\item $U_{1}\supseteq H(Y\times\left\{  1\right\}  )$ and $\overline{U}%
_{1}\cap A=\emptyset$;\smallskip

\item for $i=2,3,\cdots k$, choose $t_{i}$ so that $H(A\times\lbrack
0,t_{i}])\cap\overline{U}_{i-1}=\emptyset$ and choose $U_{i}$ containing
$H(Y\times\lbrack t_{i},1])\cup\overline{U}_{i-1}$ with $\overline{U}_{i}\cap
A=\emptyset$; and\smallskip

\item let $t_{k+1}=0$, and $U_{k+1}=Y.$\smallskip
\end{itemize}

\noindent Then $0=t_{k+1}<t_{k}<\cdots t_{2}<t_{1}=1$ and $\emptyset
=U_{0}\subseteq U_{1}\subseteq U_{2}\subseteq\cdots\subseteq U_{k}\subseteq
U_{k+1}=Y$. \smallskip

Restrict and reparametrize $H$ via the piecewise linear homeomorphism
$\lambda:[0,1]\rightarrow\lbrack0,1]$ satisfying $\lambda(0)=0$,
$\lambda(1)=1$, and $\lambda(\frac{i}{k})=t_{k-i+1}$, to get a new homotopy
$F:A\times\lbrack0,1]\rightarrow Y$ defined by $F(z,s)=H(z,\lambda(s))$.

For each $i=1,2,\ldots,k$, choose $\delta_{i}>0$ so that $B(y,\delta
_{i})\subseteq U_{i+1}$ for all $y\in U_{i}$, and let $\delta^{\prime}%
:=\min\left\{  \frac{\epsilon}{3},\delta_{i}\ |\ i=1,2,\ldots,k\right\}
$.\medskip

\noindent\textsc{Claim.}\emph{ }$F$\emph{ is a }$\delta^{\prime}$%
\emph{-}$\epsilon$\emph{-map. }\medskip

Suppose $(z,s),(z^{\prime},s^{\prime})\in F^{-1}(V)$, where
$\operatorname*{diam}V<\delta^{\prime}$. Let $y=F(z,s)$, $y^{\prime
}=F(z^{\prime},s^{\prime})$, $t=\lambda(s)$, and $t^{\prime}=\lambda
(s^{\prime})$. Choose $j\in\left\{  1,2,\ldots,k+1\right\}  $ so that $y\in
U_{j}-U_{j-1}$. Then, since $d(y,y^{\prime})<\delta^{\prime}$, we have
$y^{\prime}\in U_{j+1}-U_{j-2}$. This implies that $t_{j+1}<t<t_{j-1}$ and
$t_{j+2}<t^{\prime}<t_{j-2}$ by the choice of $t_{i}$ and $U_{i}$, so that
$\left\vert s-s^{\prime}\right\vert <\frac{3}{k}<\epsilon$. Moreover,
$d(z,z^{\prime})\leq d(z,y)+d(y,y^{\prime})+d(z^{\prime},y^{\prime
})=d(z,F(z,s))+d(y,y^{\prime})+d(F(z^{\prime},s^{\prime}),z^{\prime
})=d(z,H(z,\lambda(s)))+d(y,y^{\prime})+d(H(z^{\prime},\lambda(s^{\prime
})),z^{\prime})<\epsilon$ by definition of $H$ and the fact that
$\operatorname*{diam}V<\delta\leq\frac{\epsilon}{3}$.
\end{proof}

\end{document}